\documentclass{proc-l}   
   
\usepackage{amssymb,amscd,enumerate}   
   
\theoremstyle{plain}   
\newtheorem{theorem}{Theorem}   
   
\newtheorem{proposition}{Proposition}   
\newtheorem{corollary}{Corollary}   
   
\theoremstyle{definition}   
\newtheorem{definition}{Definition}   
   
\theoremstyle{remark}   
\newtheorem*{remark}{Remark}   
\newtheorem*{note}{Note}

\newcommand{\Lblock}[9]{\left[ \begin{array}{c|c|c}   
#1 & #2 & #3 \\   
\hline   
#4 & #5 & #6 \\   
\hline   
#7 & #8 & #9 \\   
\end{array} \right] } 
\newcommand{\homog}[2]{\hbox{\ensuremath{\mathcal{H}(2,#1,#2)}}}   
   
\newcommand{\HH}{\hbox{\ensuremath{\mathcal{H}}}}   
\newcommand{\BB}{\hbox{\ensuremath{\mathcal{B}}}}   
\newcommand{\sfi}{\hbox{\ensuremath{\mathcal{S}(\varphi)}}}   
\newcommand{\Sfi}{\hbox{\ensuremath{\mathcal{S}(\varphi)}}}   
   
\newcommand{\N}{\hbox{\ensuremath{\mathbb{N}}}}   
\newcommand{\Z}{\hbox{\ensuremath{\mathbb{Z}}}}   
   
\newcommand{\R}{\hbox{\ensuremath{\mathbb{R}}}}   
\newcommand{\C}{\hbox{\ensuremath{\mathbb{C}}}}   
\newcommand{\comment}[1]{}   
   
\DeclareMathOperator{\kker}{Ker} \DeclareMathOperator{\supp}{Supp}   
   
\numberwithin{equation}{section}   
   
\begin{document}   
   
   
   
\title
{Local Bases for Refinable Spaces}
   
   
\author[C.Cabrelli]
{Carlos Cabrelli} 
\address{
Depto.~de   
Matem\'atica \\ FCEyN\\ Univ.~de Buenos Aires\\ Cdad.~Univ., Pab.~I\\ 1428 Capital   
Federal\\ ARGENTINA\\ and CONICET, Argentina}   
\email[%
]{cabrelli@dm.uba.ar, sheinek@dm.uba.ar, umolter@dm.uba.ar} 
   
\author[S.B.Heineken]{Sigrid B.~Heineken}   
   
\author[U.M.Molter]{Ursula~M.~Molter}   
   
\thanks{The research of   
the authors is partially supported by Grants:   
CONICET, PIP456/98, and UBACyT X058 and X108. The authors also acknowledge partial support  from the Guggenheim Foundation and the Fulbright Commission during the period in which part of this research was performed.}   
   
   
\keywords{Homogeneous functions, shift-invariant spaces, accuracy, refinable functions}  
\subjclass{Primary:39A10, 42C40, 41A15}   
   
\date{\today}   
   
   
\begin{abstract}   
We provide a new   
representation of a refinable shift invariant space with a compactly
supported generator, in terms of    
functions with a special property of homogeneity. In particular
these functions include all the homogeneous polynomials
that are reproducible by the generator, what links this representation
to the accuracy of the space. 
We completely characterize the class of homogeneous functions in
the space and show that they can reproduce the generator. 
As a result we conclude that the homogeneous functions can
be constructed from the vectors associated to
the spectrum of the scale matrix (a   
finite square matrix with entries from the mask of the   
generator).  Furthermore, we prove
that the kernel of the transition operator has
the same dimension than the kernel of this finite matrix.
This relation provides an easy test for the linear independence
of the integer translates of the generator.
This could be potentially   
useful in applications to approximation theory, wavelet theory and   
sampling.    
\end{abstract}   
   
\maketitle   
   
   
\section{Introduction}

A function  $\varphi: \R \longrightarrow \C$ is called {\em refinable} if    
it satisfies the equation:   
\begin{equation}\label{dilation-equation}   
\varphi(x) = \sum_{k=0}^N c_k \varphi(2x - k),    
\end{equation}   
for some complex scalars $c_0,...,c_N$.   
The scalars ${c_k}$ are the {\em mask} of the refinable function.   
We consider the case in which $\varphi$ is compactly supported.   
Define the {\em Shift Invariant Space} (SIS) generated by $\varphi$   
as:    
$$ {\mathcal S}(\varphi) = \{f: \R \longrightarrow \C: f(x) =   
\sum_{k\in \Z} y_k \varphi(x+k), y_k \in \C \}.   
$$   
A refinable SIS is a SIS with a refinable generator.   
Refinable SIS and refinable generators  have been studied extensively, since   
they are very important in Approximation Theory and Wavelet Theory.   
   
Many properties of  $\varphi$  can be obtained imposing conditions   
 on the mask.
One fundamental question is when the space $\sfi$ contains polynomials   
 and of which degree. The {\em accuracy} of $\varphi$ is the maximum integer $n$ such that   
all the polynomials of degree less or equal than $n-1$ are   
contained in $\sfi $.   
   
The accuracy is related to the approximation order of   
$\sfi$, (\cite{Jia95}, \cite{deB90} and references therein), 
and with the zero moments 
 and the  smoothness of the associated wavelet when $\varphi$ generates   
 a Multiresolution Analysis \cite{Mey92}.   
There are many well known equivalent conditions for accuracy.   
The one that interests us here is the following, 
\cite{Dau92}, \cite{CHM98}:   
\begin{proposition} \label{accuracy}  
Let $\varphi$ be a compactly supported function satisfying   
\eqref{dilation-equation}. Then 
$\varphi$ has accuracy $n$ if and only if  
$\{1,2^{-1},...,2^{-(n-1))}\}$ are   
 eigenvalues of the $(N+1)\times (N+1)$ matrix $T$ defined by  
$T = \{c_{2i-j}\}_{i,j=0,...,N}$, (the {\em scale matrix})  and   
there exist polynomials $p_0,...,p_{n-1}$ with degree$(p_i)=i$   
such that the vector $v_i= \{p_i(k)\}_{k=0,...,N}$     
is a left eigenvector of $T$ corresponding to the eigenvalue   
$2^{-i}$.    
\end{proposition}   
Here, and always throughout the paper, we assume $c_t = 0$, if $t \not= 0,\dots, N.$ 
One interesting property is that if $\varphi$ has accuracy $n$,   
then for $s=0,1,...,n-1$ it is true that   
$x^s = \sum_{k\in\Z} p_s(k) \varphi(x-k),   
$
where $p_s$ is the polynomial that provides the eigenvector   
for the eigenvalue $2^{-s}$.   
So, if the polynomial $x^s$ is in $\Sfi$, then the left eigenvector of $T$, 
corresponding to the eigenvalue $2^{-s}$, provides the coefficients needed
to write $x^s$ as a linear combination of the translates of $\varphi$.

A {\em local basis} for $\Sfi$ is a set of functions in $\Sfi$ whose restriction
to the $[0,1]$ form a basis for the space of
all the functions in $\Sfi$, restricted to the $[0,1]$.

When $\varphi$ is the B-spline of order $m$ (so $N=m$), then all
polynomials of degree less or equal than $m-1$ are in $\Sfi$. 
Moreover, the set $\{1, x, x^2, \dots, x^{m-1}\}$ is a local basis
for $\Sfi$, and the spectrum of $T$ consists exactly of
$\{1, 2^{-1}, \dots, 2^{-(m-1)}\}$.

Now, if $\varphi$ is not a B-spline then $T$ could have some
eigenvalue $\lambda$ different from a power of $1/2$.
If powers of $1/2$ are associated to homogeneous polynomials, which
functions in $\Sfi$ are associated to an arbitrary eigenvalue $\lambda$? 
Will the functions, associated to all the eigenvalues, provide also
local bases of $\Sfi$, or equivalently, will these functions reproduce
the generator $\varphi$? 

In the case that $\lambda$ is a simple eigenvalue, Blu and Unser
\cite{BU02} and later Zhou \cite{Zho02} showed 
that $\lambda$ is associated to what they call a $2$-scale $\lambda$-homogeneous
function, that is a function in the SIS that satisfies the relation
$h(x) = \lambda h(2x)$.

However, to obtain a complete representation of the space it is
necessary to consider the whole spectrum of $T$. 
This motivates the study of the spectral properties of $T$ for a {\em general}
refinable $\varphi$.
This is achieved in this paper: we are able to completely characterize the
SIS in terms of functions associated to the spectrum of $T$.
We prove that these functions provide a local basis of $\Sfi$ (c.f. Theorem~\ref{main}).
The advantage of this local basis is that it contains all possible
homogeneous polynomials in the space, and those functions in the
basis which are not polynomials, still preserve some
kind of homogeneity. Furthermore this basis can be easily
obtained from the spectrum of the finite matrix $T$. We are also able to 
prove that this matrix is necessarily invertible if the translates of $\varphi$
are linearly independent.
\begin{definition} \label{defi-hom}
Let $\lambda \in \C,  \lambda \not= 0$ and $r\geq 1$ an integer.   
A function $h$ is $(2,\lambda,r)$-{\em homogeneous} if it satisfies   
the following equation:   
\begin{equation} \label{homog-eq}  
\sum_{k=0}^r {r \choose k}   (- \lambda)^{r-k} h(2^{-k} x) = 0   
\;\; \textrm{a.e.}.   
\end{equation}  
$r$ is called the order of homogeneity,   
and $\lambda$ the degree.
\end{definition}

If $\HH \subset \Sfi$ is the span of all the
$(2,\lambda,r)$-homogeneous functions in $\Sfi$, for any $\lambda \in \C$ and
any positive integer $r$, we show that under the hypothesis of linear independence
of the translates of $\varphi$, $\text{dim} (\HH) = N+1$, and that there is a
basis of $\HH$, corresponding to the spectrum of $T$. 
More precisely, given a basis $\BB = \{v_0, \dots, v_{N}\}$ of $\C^{N+1}$ that
yields the Jordan form of $T$ we associate to each vector $v \in \BB$ a unique
$(2,\lambda,r)$-homogeneous function in $S(\varphi)$, where $\lambda$ and $r$
satisfy   $v(T-\lambda I)^r = 0$.

The first $N$ of these functions are a local basis of functions 
in $S(\varphi)$ restricted to $[0,1]$.   
This allows to reconstruct   
the generator $\varphi$ from the homogeneous functions,
and gives a new representation for the functions in $S(\varphi)$.   

Furthermore, we show that to 
each non-zero vector in the kernel of $T$, there corresponds a non-trivial linear 
combination of the integer translates of $\varphi$ yielding the zero function. 
   
\section{Notation}  \label{notation} 
Let $\varphi: \R \longrightarrow \C$ be a function supported in $[0,N]$ 
satisfying \eqref{dilation-equation}.  
   
We will often use an infinite column vector associated to $\varphi$,   
namely    
\begin{equation}   
\left(\phi(x)\right)^t = \left[ \dots,\ \varphi(x-1),\ \varphi(x), \  
\varphi(x+1),\ \dots\ \right]. \label{vfi}   
\end{equation}   
Let $\ell(\Z)$ be the space of all the sequences defined in $\Z$.   
We say that the integer translates of $\varphi$ are {\em globally linearly independent}, 
or {\em linearly independent} if  
$\sum_{k\in \Z} \alpha_k \varphi(\cdot - k) \equiv 0 \quad \Longrightarrow \quad \alpha_k = 0\   
\forall k$, 
for any sequence $\alpha \in \ell(\Z)$.  
 
The {\em subdivision operator} associated to the mask $c_k$ is   
the operator   
\begin{equation} 
\label{subdiv}   
 S_c :\ell(\Z)\rightarrow \ell(\Z)   
\quad \text{defined by} \quad   
S_c(\alpha)_j = \sum_{i \in \Z} \alpha_i c_{2i-j}.   
\end{equation} 

\begin{note} The subdivision operator is sometimes defined in a different but equivalent  
way as $\tilde{S}_c(\alpha)_j = \sum_{i\in \Z} \alpha_i c_{j-2i}$. If $h: \ell(\Z)
\rightarrow \ell(\Z)$ is the 
operator $h(\alpha)_k = \alpha_{-k}$, then $h\tilde{S}_c h = S_c$ and
therefore  $S_c$ and $\tilde{S}_c$ 
share most of the properties. 
For a nice account of properties of the subdivision operator
see \cite{BJ02}. \end{note} 
   
If $L = L_{\varphi}$ is the double infinite matrix   
$L = \left[c_{2i-j}\right]_{i,j \in \Z}$,   
then the refinement equation can be written as   
$\phi(x) = L \phi(2x)$.   
 
Using the matrix $L$, the subdivision operator \eqref{subdiv} can be recast as: 
$S_c \alpha = \alpha L,\ \alpha \in \ell(\Z),$ 
where $\alpha$ on the right hand side of the equation is thought as an infinite row vector. 
Note that the scaling matrix $T$ defined in Proposition~\ref{accuracy} 
is a finite submatrix of $L$. 
We will consider in our analysis the matrices $M, T_0, T_1$, that are  
submatrices of $T$ and are defined as: 
$M=[c_{2i-j}]_{i,j=1,...,N-1},\ T_0=[c_{2i-j}]_{i,j=0,...,N-1},\ 
T_1=[c_{2i-j}]_{i,j=1,...,N}$. That is, 
\begin{equation} \label{Tone}   
T = 
\left[ \begin{array}{c|c|c}   
c_0 & 0 & 0 \\   
\hline    
\vdots & M & \vdots\\   
\hline    
0 & 0 & c_N    
\end{array} \right] .   
\end{equation}   
 Note that $c_0$ and $c_N$ must be non-zero, 
since $\supp(\varphi) = [0,N]$.

Now, if $Y \in \ell(\Z)$,   
define $Y^0$ and $Y^M$ as the restriction of $Y$ to   
the indexes $\{0, \dots, N\}$,   
and   
$\{1, \dots, N-1\}$, respectively, i.e.,    
$$Y^0 = (Y_0, \dots, Y_{N}), \quad Y^M = (Y_1, \dots, Y_{N-1}). $$
\begin{note}   
Throughout this paper, $(L - \lambda I)$ is considered as an operator on   
$\ell(\Z)$, defined by left-multiplication, (i.e. $ Y \longmapsto   
Y(L-\lambda I)$, where $Y$ is a double infinite row vector).   
$I$ is the identity operator acting on $\ell(\Z)$. By an abuse of   
notation, we will use   
the notation $I$ for all identity operators, without   
distinguishing the space they are acting on. 
\end{note}   
   
\section{The point spectrum of $L$}   
  
The following proposition, will show, how the spectral properties of $L$    
are related to those of $T$.  The case $r=1$ has been studied earlier by \cite{CHM00}, 
\cite{JRZ98}, \cite{Zho00,Zho02}.  
\begin{proposition} \label{I}   
Let $\lambda \in \C.$   
\begin{enumerate}   
\item \label{p1}   
Let $Y \in \ell(\Z)$ and $r \in \N$, $r \geq 1$. 
If $Y \in \kker(L - \lambda I)^r$, then $   
Y^0 \in \kker(T - \lambda I)^r.$    
Moreover, if    
$\lambda \not= 0$, $Y \not= 0$ and $Y \in \kker(L - \lambda I)^r$, then   
$Y^0 \not= 0$.   
\item If $v \in \kker(T - \lambda I)^r$ and $\lambda \not= 0$,   
then there exists an extension $Y_v \in \ell(\Z)$ of $v$,   
(i.e. $Y_v^0 = v$) such that $Y_v \in \kker(L - \lambda I)^r$. \label{p2}   
\end{enumerate}   
\end{proposition}   
\begin{proof}   
The matrix $L$ (and therefore $L-\lambda I$) can be decomposed in blocks as   
\begin{equation} \label{L-blocks}   
L = \Lblock{R}{0}{0}{P}{T}{Q}{0}{0}{S},\quad 
(L - \lambda I) = \Lblock{R - \lambda I}{0}{0}{P}{T - \lambda I}{Q}%
{0}{0}{S -\lambda I},   
\end{equation}   
where we decompose $\Z$ as   
$\Z =  A^{-}  \cup  A^0    \cup  A^{+},$    
with   
$A^{-} = \Z  \cap (-\infty,-1]$, $A^0 = \Z  \cap  [0,N]$ and   
$A^{+} = \Z  \cap [N+1, +\infty)$,   
and   
\begin{equation*}   
R = L|_{A^-\times A^-} \quad P = L|_{A^0\times A^-} \quad T =   
L|_{A^0\times A^0} \quad Q = L|_{A^0\times A^+} \quad S =   
L|_{A^+\times A^+}.   
\end{equation*}   
This block form of the matrix, is closed under multiplication. So if $r \geq 1,   
r \in \N$   
\begin{equation} \label{PQ}   
L^r = \Lblock{R^r}{0}{0}{P_r}{T^r}{Q_r}{0}{0}{S^r}, \quad \text{and} \quad   
(L - \lambda I)^r = \Lblock{(R - \lambda I)^r}{0}{0}%
{P_r^{\lambda}}{(T - \lambda I)^r}{Q_r^{\lambda}}%
{0}{0}{(S -\lambda I)^r} ,   
\end{equation}   
where   
$P_r = \sum_{k=0}^{r-1} T^k P R^{r-k-1}$ and $
Q_r = \sum_{k=0}^{r-1} T^k Q S^{r-k-1},$   
and $P_r^{\lambda}$ and $Q_r^{\lambda}$ are analogous, with   
the obvious changes.    
   
Note that the matrix $S$ is upper triangular,    
with diagonal $(0, 0, 0, \dots)$, and hence    
$(S - \lambda I)^r$ is upper triangular, with diagonal   
$((-\lambda)^r, (-\lambda)^r, (-\lambda^r), \dots)$.   
   
Analogously, we observe that $R$ is lower triangular with zeroes in the    
main diagonal, so  $(R -\lambda I)^r$ is   
lower triangular with diagonal    
$((-\lambda)^r,(-\lambda)^r,(-\lambda)^r,\dots)$.

If $Y = (Y^-,Y^0,Y^+)$, then   
$Y(L - \lambda I)^r$ can be written as
\begin{equation}    
(Y^-(R - \lambda I)^r + Y^0 P_r^{\lambda},    
Y^0(T - \lambda I)^r,   
Y^0 Q_r^{\lambda} + Y^+(S - \lambda I)^r). 
\label{ydescomp}   
\end{equation}   
So if $Y \in \kker(L - \lambda I)^r$, then $Y^0 \in \kker(T - \lambda I)^r$.   
   
We now want to show that if $Y \in \kker(L - \lambda I)^r$,    
$\lambda \not= 0$, $Y \not= 0$, then $Y^0 \not= 0$.    
   
For this, let $k_0 \in \Z$    
be such that $Y_{k_0} \not= 0$. If $0 \leq k_0 \leq N$, we are   
done. Assume that $k_0 > N$. Then, since   
$Y(L-\lambda I)^r = 0$, in particular, the $k_0$ element of this   
product is $0$.    
But since $\lambda \not=0$, $(S - \lambda I)^r$ is   
upper triangular with $(-\lambda)^r$ in the diagonal,    
therefore the only nonzero elements of   
 column $k_0$ of $(L - \lambda I)^r$ are between $0$ and $k_0$.   
Hence there has to be a $k_1, 0 \leq k_1 < k_0$ such that   
$Y_{k_1} \not= 0$. Again, if $0 \leq k_1 \leq N$ we are   
done, otherwise we repeat the argument until $k_j$ is in the    
desired interval.   
If $k_0 < 0$, the argument works in the same way, reversing   
the role of $(S-\lambda I)^r$ and $(R - \lambda I)^r$.   
   
For the proof of part~\ref{p2}, assume that $v \in \C^{N+1}$,    
$v \in \kker(T - \lambda I)^r$.    
We want to find an infinite  vector $Y \in \ell (\Z)$, such that    
$Y^0 = v$ and $Y \in \kker(L - \lambda I)^r$.   
From equation \eqref{ydescomp} we know that if $Y \in \ell (\Z)$,   
\begin{equation*}   
\left[Y (L - \lambda I)^r \right]^+  =  
 Y^0 Q_r^{\lambda} + Y^+(S-\lambda I)^r,\quad   
\left[Y (L - \lambda I)^r \right]^-  =  Y^0 P_r^{\lambda} + Y^-(R-\lambda I)^r.   
\end{equation*}   
Therefore, if $ Y \in \kker(L - \lambda I)^r$, and $Y^0 = v$, then $Y^+$ and $Y^-$ have 
to satisfy   
\begin{equation}\label{4*}   
Y^+(S-\lambda I)^r = -v Q_r^{\lambda}  
\quad \text{and} \quad 
Y^-(R-\lambda I)^r = -v P_r^{\lambda} .   
\end{equation}   
Using again that $(S - \lambda I)^r$ and $(R - \lambda I)^r$, are triangular, if   
$\lambda \not= 0$, there are unique solutions for $Y^+$ and $Y^-$ and they can be   
obtained recursively.   
\end{proof}   
The last proposition tells us that the elements of the spectrum of   
$T$ are intimately related to those of $L$. But by the special form   
of $T$ (see equation \eqref{Tone}), we can actually use the $(N-1) \times
(N-1)$ matrix $M$ to obtain the spectrum of $T$, as the following proposition
shows:
\begin{proposition} \label{IIbis}   
Let $\lambda \not=0 \in \C$.    
\begin{enumerate}   
\item Let $v^0 = (v_0,\dots, v_N) \in \C^{N+1}$ and $r \in \N$, $r \geq 1$.    
Then,   
if $v^0 \in \kker(T - \lambda I)^r$ then $v^M = (v_1, \dots, v_{N-1})   
 \in \kker(M - \lambda I)^r$.    
Moreover, if $\lambda \not= c_0$, $\lambda \not= c_N$,   
$v^0 \in \kker(T - \lambda I)^r$, and   
$v^0 \not= 0$, then   
$v^M \not= 0$.   
\item if $v^M = (v_1,\dots,v_{N-1}) \in \kker(M - \lambda I)^r$ and    
$\lambda \not= c_0$ and $\lambda \not= c_N$,   
then there exists an extension $v^0 \in \C^{N+1}$ of $v$,   
such that $v^0 \in \kker(T - \lambda I)^r$.   
\end{enumerate}   
\end{proposition}   
\begin{proof}   
The proof is immediate by noting the special block-form of $T$. 
\end{proof}   

\section{The kernel of $L$}
   
The case $\lambda = 0$ could not be handled with the   
methods of Proposition~\ref{I},   
since the matrices $R$ and $S$ in \eqref{L-blocks} have zeros in   
the main diagonal. Instead,    
 we need some results from   
the theory of difference equations which we present below \cite{Hen62}.   
   
\subsection{Difference Equations}   
   
Consider the linear difference equation with constant coefficients of order $r$   
\begin{equation}   
u_0y_n + u_1 y_{n+1} + \dots + u_r y_{n+r} = 0  \quad y = \{y_n\}_{n \in \Z},   
\label{diff}   
\end{equation}   
where $u_k \in \C, u_0 \not=0, u_r \not= 0$ with characteristic polynomial   
$P(x) = \sum_{k=0}^r u_{k} x^{k}.   
$
   
A {\em solution} to the equation~\eqref{diff} is  a sequence $Y$ in $\ell(\Z)$, that 
satisfies \eqref{diff}   for all $k \in \Z$. A vector $y=(y_0, \dots, y_m)$ with 
$m \geq r+1$ is a {\em finite solution} of \eqref{diff}, if it satisfies 
\eqref{diff} for $n=0$ to $n=m-r-1$. 
 
The {\em space of solutions}  $S \subset \ell(\Z)$, has dimension $r$, and a basis   
of this space (the fundamental basis) can be written in the following way:   

Let $h \geq 1$ be an integer, $d_1,\dots,d_h$ arbitrary non-zero complex numbers 
with $d_i \not= d_j$ if $i \not= j$. Let $r_1,\dots,r_h$ be positive integers. To each pair 
$(d_i,r_i)$, $i=1,\dots,h$ we will associate a sequence $a_i = \{a_{ik}\}_{k\in \Z}$ 
defined as follows: 
$\text{Set} \quad r = r_1+\dots+r_h \quad \text{and}\quad r_0 = 0.   
$\ Let $0 \leq i \leq r-1$ and $s = s(i), j = j(i)$ be the unique integers that   
satisfy   
$
r_0+\dots+r_{s-1} \leq i < r_1+\dots+r_s,\ j(i) = i - \sum_{k=0}^{s(i)-1}r_k$.   
\begin{equation}\label{aik}   
\text{Define} \quad a_{ik} =    
\begin{cases}   
\text{sg}(k)\frac{|k|!}{(|k|-j(i))!}\, d^k_{s(i)} & \text{for}\ |k|  \geq j(i)\\   
0 & |k| < j(i).   
\end{cases} \quad i=0,\dots,r-1, \quad k \in \Z,
\end{equation}   
where $\text{sg}(k)$ is the sign of $k$.

So, if $P$ is the characteristic polynomial associated to~\eqref{diff}, 
consider the pairs $\{(d_i,r_i): \text{where}\ d_i\ \- \text{is a root of $P$ and $r_i$  
its multiplicity} \}$. 
The sequences $\{a_{ik}\}_{k \in \Z}$, $i = 0,\dots,r-1$ form a basis   
of $S$, the subspace of $\ell(\Z)$ of the solutions to \eqref{diff}.   
 
It is also known from the theory of difference equations, that every solution is 
determined unequivocally by any $r$ consecutive elements of it. Hence, if $y$ is a solution 
such that $r$ consecutive elements are $0$, then $y$ is the zero solution. 
 
We will now associate to the pairs $\{(d_i,r_i): i = 1,\dots,h\}$ the $r\times r$ 
matrix $A = [a_{ij}]_{i,j = 0, \dots, r-1}$. Then (cf. Henrici \cite{Hen62}, pg.~214)  
\begin{equation} \label{deter}   
\det(A) = \prod_{1 \leq l < s \leq h} (d_l - d_s)^{r_l + r_s}   
\prod_{i=1}^h (r_i - 1)!! ,   
\end{equation}   
where $0!! = 1$ and $k!! = k!(k-1)!\dots 1!$.   
Since $d_i \not= d_j$ for $i\not=j$, $\det(A) \not=0$ and $A$ is invertible.  
 
Let us now 
consider a system of $k$ linear difference equations with constant coefficients of order $r$.   
\begin{equation} \label{sys-diff}   
u_{i0}y_n + \dots + u_{ir} y_{n+r} = 0 
, \quad i = 1, \dots k, \quad n \in \Z,   
\end{equation}   
and let $P_i$ be the characteristic polynomial of equation $i$, 
$P_i(x) = \sum_{j=0}^r u_{ij}x^j$. Define $
\mathcal{D} = \cup_{i=1}^k \{d : P_i(d) = 0\} = \{d_1,\dots,d_s\}$, 
and for each $d \in \mathcal{D}$ define 
$r_d = max\{r_i: r_i \text{ is the multiplicity of $d$ in $P_i$}\}$. 
Note that $r_d \geq 1\ \forall d \in \mathcal{D}$. We then have the pairs 
$(d_i,r_{d_i}) = (d_i,r_i)$. Define the index of the system 
to be  $t = \sum_{d\in \mathcal{D}} r_d \leq kr$. 
Let $\ell$ be the degree of the maximum common divisor $p$, of 
$\{P_a,\dots,P_k\}$. Hence, $P_i(x) = p(x) \tilde{P}_i(x)$, with degree $\tilde{P}_i = r-\ell$. 
(Note that $\ell$ could be $0$). 
With the above notation, we have the following proposition:
\begin{proposition}\label{pdiff}   
The space $S_k$ of solutions to the system~\eqref{sys-diff} 
has dimension $\ell$, where $\ell$ is the degree of the maximum common divisor of 
the characteristic polynomials. Moreover, if $t$ is the index of the system~\eqref{%
sys-diff}, and
$z$ is a vector of length $t$ that satisfies~\eqref{sys-diff}, then it can be 
extended to a sequence $y_{z} = \{y_j\}_{j\in \Z}$ solution  of \eqref{sys-diff} and such that 
$y_j = z_j, \,\, j=1,\dots,t$. 
\end{proposition}   
\begin{proof}   
Let $p$ be the maximum common divisor of $P_1,\dots,P_k$, and let $\ell = \deg(p)$. 
It is clear, that $\dim(S_k) \geq \ell$.  
 
For the other inequality,  
consider the $t \times t$ matrix $A = [a_{ij}]_{i,j = 0, \dots, t-1}$,   
with $a_{ij}$ defined in \eqref{aik} for the pairs 
$\{(d_i,r_i)\}$ defined above and $t$ being the index of the 
system~\eqref{sys-diff}.  
Since $d_i \not= d_j$ for $i\not=j$, $\det(A) \not=0$ by  
\eqref{deter} and $A$ is invertible.  
 
Assume now that $y \in S_k$, then $y$ is a solution to all $k$ difference equations, 
hence there exist 
$\alpha^1, \dots, \alpha^k$ vectors of length $r$, such that 
\begin{equation} 
A_i \alpha^i = \left[
y_0, \dots, y_{t-1}\right]^t  
\quad 1\leq i \leq k, 
\end{equation} 
where $A^i$ is an $t \times r$ matrix whose columns are a fundamental system for 
equation $i$. Note that $A^i$ is a sub-matrix of $A$, whose columns correspond 
to some columns $\{i_1,\dots,i_r\}$ of $A$. 
 
Let now $\tilde{\alpha}^i$ be vectors of length $t$, such that $\tilde{\alpha}^i_h = 0$ 
whenever $h \not\in \{i_1,\dots,i_r\}$ and $\tilde{\alpha}^i_{i_s} = \alpha_s$, 
$s = 1, \dots, r$. Then we have for $i, j = 1,\dots,k$ 
\begin{equation} 
A \tilde{\alpha}^i = A_i \alpha^i = \left[ 
y_0, \dots, y_{t-1}\right]^t = A_j \alpha^j =  A \tilde{\alpha}^j, 
\end{equation} 
and hence  
$A (\tilde{\alpha}^i - \tilde{\alpha}^j) = 0$, for all $i \not= j$ 
and therefore, by the invertibility of $A$, $\tilde{\alpha}^i = \tilde{\alpha}^j$, for all $i \not= j$. Therefore the only non-zero elements of $\alpha_i$ can be those corresponding 
to the columns associated to the roots of $p$. Hence $y$ is a linear combination 
of $\ell$ columns, and therefore $\dim(S_k) \leq \ell$. 
\end{proof}   
 
By noting that for the previous proof, we only used the 
first $t$ coordinates of the infinite sequences, we have the following immediate Corollary.  
\begin{corollary} 
\label{pp2} If $z$ is a vector of length $t$ that satisfies~\eqref{sys-diff}, then it can be 
extended to a sequence $y_{z} = \{y_j\}_{j\in \Z}$ solution  of \eqref{sys-diff} and such that 
$y_j = z_j, \,\, j=1,\dots,t$. 
\end{corollary}   
\subsection{The $\kker(L)$}   
   
We can now return to our double infinite matrix $L$ and look at the special case $\lambda = 0$.    
As it turns   
out, the kernel of $L$ is characterized by the vectors in the   
kernel of $M$. Since $c_0$ and $c_N$ are non-zero, the matrices   
$T$ and $M$ have kernels of the same dimension.   
Moreover, we have the following Proposition: 
\begin{proposition} \label{II} 
Consider the polynomials $p_e$ and $p_o$  
of degree $q = \frac{N-1}{2}$ (we assume $N$ to be odd)   
$ p_e(x) = c_0 + c_2 x + \dots + c_{2q} x^q, \ 
p_o(x) = c_1 + c_3 x + \dots + c_{2q+1} x^q.$ 
Then $\dim(\kker(L)) = \dim(\kker(M)) = \text{\rm degree}(p)$, 
where $p$ is the  
maximum common divisor of the polynomials $p_e$ and $p_o$. 
In particular, if $\dim(\kker(M)) >0$, $p_e$ and $p_o$ 
have a common root. Furthermore 
\begin{enumerate} 
\item \label{one} 
For every $Y \in \kker(L)$, $Y \not= 0$, we have $Y^M \not=0$ 
and $Y^M \in \kker(M)$. 
\item \label{two} 
Conversely, for each $v \in \kker(M), v \not= 0$, we have 
$Y_v \not= 0$ and $Y_v \in \kker(L)$.   
\end{enumerate} 
\end{proposition}   
\begin{proof}   
Let us observe first, that $Y \in \ell(\Z)$ is in the Kernel of $L$, if and 
only if $Y$ satisfies the system of difference equations: 
\begin{equation} \label{eod} 
\left\{ 
\begin{array}{ll} 
c_0v_n + c_2 v_{n+1} + \dots + c_{2q} v_{n+q} & = 0  \\   
c_1v_n + c_3 v_{n+1} + \dots + c_{2q+1} v_{n+q} & = 0. 
\end{array} 
\right. 
\end{equation} 
Therefore, by Proposition~\ref{pdiff}, $\kker(L)$ is the subspace generated 
by the fundamental solutions associated to the roots of $p$, the maximum 
common divisor of $p_o$ and $p_e$. This shows that $\dim(\kker(L)) = \text{degree}(p)$. 
 
On the other side, if $Y \in \kker(L)$, since $(YL)^M = Y^M M$ we conclude 
that $Y^M \in \kker(M)$, and if $Y^M$ is the zero vector, then the solution 
$Y$ of \eqref{eod} has $N-1$ consecutive zeros, so $Y = 0$. 
Hence, if $Y \not= 0$, then $Y^M \not= 0$, which proves~(\ref{one}). 
    
To see that if $v=(v_1,\dots,v_{N-1})$ satisfies $vM = 0$, then $v$ can be extended, just 
note that the   
sequence $v_1, \dots, v_{N-1}$ must satisfy the difference equations system  
of order $q = \frac{N-1}{2}$ (we assumed $N$ to be odd) given by \eqref{eod}. 
Since the index $t$ of the system~\ref{eod} satisfies $t \leq 2q = N-1$,  
and $v$ is a non-trivial common solution of length $N-1$, by Corollary~\ref{pp2},  
this solution can   
be extended in such a way that the extension satisfies both difference equations. 
This proves~(\ref{two}). 
 
From~(\ref{one})~and~(\ref{two}) it is immediate that $\dim(\kker(L)) = \dim(\kker(M))$.   
\end{proof}   

\begin{note} The fact that $\dim(\kker(M)) >0$ implies that $p_e$ and $p_o$ 
have a common root was proved under some minor technical conditions by Meyer \cite{Mey91}. Related results can also be found in \cite{JW93}.
\end{note}

\subsection{Invertibility of $L$}   
   
Propositions~\ref{I} and \ref{IIbis}, relate the spectral properties   
of the matrix $M$ to the ones of the operator $L$. The next proposition   
shows a necessary condition for the independence of the integer translates of   
the function $\varphi$, in terms of the matrix $M$.   
\begin{proposition} \label{prop6}   
With the above notation, consider the following properties   
\begin{enumerate}   
\item $\{\varphi(\cdot - k)\}_{k\in \Z}$ are globally linearly independent, \label{iI}   
\item The operator $L:\ell(\Z) \longrightarrow \ell(\Z),   
Y \longmapsto YL$ is one-to-one, \label{iII}   
\item The matrix $M$ is invertible. \label{iIII}   
\end{enumerate}   
Then (\ref{iII} $\Longleftrightarrow$ \ref{iIII}) and (\ref{iI}    
$\Longrightarrow$ \ref{iII}).   
\end{proposition}   
\begin{proof}   
{\bf (\ref{iI} $\Longrightarrow$ \ref{iII})}   
Assume $YL = 0$. Define $F(x) = Y \phi(x)$. Then we have   
$F(x) = Y\phi(x) = YL\phi(2x) = 0.$
Now, $Y\phi(x) = 0  \Longrightarrow Y = 0$, therefore $\kker(L) = \{0\}$.   
   
\noindent   
{\bf (\ref{iII} $\Longleftrightarrow$ \ref{iIII})}   
is a consequence of Proposition~\ref{II}.    
\end{proof}   
   
\noindent   
{\bf Note:} We do not know if either   
 (\ref{iII}) or   
(\ref{iIII}) implies (\ref{iI}).    
   
\section{Homogeneous functions}   
   
Assume now that $Y \in \kker(L-\lambda I)^r$, and define the function $
h \in \Sfi$ as
$h(x) = Y\phi(x)$. So $h$ is $(2, \lambda,r)$ homogeneous (c.f.~\eqref{homog-eq}),
since it satisfies:   
\begin{align*}   
0 &= Y (L - \lambda I)^r\phi(x) =    
Y \left(\sum_{k=0}^r \begin{pmatrix}r\\k   
\end{pmatrix} (-\lambda)^{k} L^{r-k}\right) \phi(x) \\   
 &= Y \left(\sum_{k=0}^r \begin{pmatrix}r\\k\end{pmatrix}   
 (-\lambda)^{k} \phi(\frac{x}{2^{r-k}})\right)   
 = \sum_{k=0}^r \begin{pmatrix}r\\k\end{pmatrix}   
 (-\lambda)^{k} h(\frac{x}{2^{r-k}}). 
\end{align*}   
We will denote by   
\homog{\lambda}{r}, the space of all $(2,\lambda,r)$ homogeneous   
functions.   
\begin{remark}[1]   As pointed out by the referee, a $(2,\lambda,r)$-homogeneous function, is a particular case of a poly-scale refinable distribution. The concept of poly-scale refinable distribution is a generalization of refinability. See for instance \cite{DD02, Sun05}.
\end{remark}
\begin{remark}[2]   
Note that if    
$h \in \homog{\lambda}{r}$  then   
$h \in \homog{\lambda}{s}$ for every $s \geq r$.  
Therefore the ``order of homogeneity'' will be defined by    
$\min\{s: h \in \homog{\lambda}{s}\}$.   
\end{remark}    
\begin{remark}[3] If $h$ is homogeneous (of any order) and $\lambda \not= 1$, then $h(0) = 0$.    
The values of any homogeneous function of order $r$   
in $(0,+\infty)$, are completely determined   
by its values on any interval of the type    
$\left[\frac{1}{2^{k+r}}, \frac{1}{2^k}\right)$,   
$k \in \Z$.    
(Analogously, the values on $(-\infty,0)$, are obtained from the values   
in any interval of the type   
$\left(-\frac{1}{2^{k}}, -\frac{1}{2^{k+r}}\right]$).   
\end{remark}   
\begin{remark}[4]   
In the case of order of homogeneity $1$, (e.g. $r=1$), $h$ is a $2$-scale   
homogeneous function as described  in \cite{Zho02}.   
\end{remark}   
\begin{proposition}   
Assume $\{\varphi(\cdot-k)\}$ are linearly independent. Let   
$\phi$ be as in \eqref{vfi}. If $g_1, \dots, g_n \in \sfi$,   
$g_i = Y^i \phi$, then   
$\{g_1, \dots, g_n\}$ are linearly independent functions   
if and only if 
$ \{Y^1, \dots, Y^n\}$ are linearly independent in $\ell(\Z)$.   
\end{proposition}   
\begin{proof}   
We observe that   
$\sum_{i=1}^n \alpha_i g_i = \sum_{i=1}^n \alpha_i\left(Y^i \phi\right) =   
\left(\sum_{i=1}^n \alpha_i Y^i \right)\phi.$    
This equation, together with the linear independence of the translates   
of $\varphi$, tells us that   
$\sum_i \alpha_i g_i \equiv 0$ if and only if    
$\left(\sum_i \alpha_i Y^i \right)= 0$, which proves the desired result.   
\end{proof}   
\begin{theorem}\label{tone}   
Assume that $\{\varphi(\cdot-k)\}$ are linearly independent. If $h = Y\phi$,   
($ h \in \sfi$), and $h \in \homog{\lambda}{r}$   
then $v_h = Y^0 \in \kker(T-\lambda I)^r$. Reciprocally, if   
$v \in \kker(T-\lambda I)^r$, then the function $h = Y_v \phi$ is   
in \homog{\lambda}{r}. (Here $Y_v$ is the unique extension   
of $v$ to a vector in $\kker(L - \lambda I)^r$ by Prop.~\ref{I}.)   
\end{theorem}   
\begin{proof}   
For the first claim, note that   
\begin{equation*}   
0  = \sum_{k=0}^r \begin{pmatrix}r\\k\end{pmatrix} (-\lambda)^{r-k} h(2^{-k} x)   
 = \sum_{k=0}^r \begin{pmatrix}r\\k\end{pmatrix} (-\lambda)^{r-k}%
YL^k \phi(x)\\   
 = Y \left(L - \lambda I\right)^r \phi(x).   
\end{equation*}   
Then, $ Y \left(L - \lambda I\right)^r=0$ and by 
 Proposition~\ref{I}, $v_h \in \kker(T- \lambda I)^r$.   
   
For the converse first observe that if $v=0$ the result is trivial. Assume 
$v\not= 0$ and $v \in \kker(T- \lambda I)^r$, then by   
Proposition~\ref{prop6} $\lambda \not= 0$. Hence (by Prop.~\ref{I}) there is a unique extension 
 $Y_v \in \kker(L- \lambda I)^r$, so   
$h = Y_v \phi$ is in \homog{\lambda}{r}.   
\end{proof}    
   
\subsection{Local basis of homogeneous functions}   
\begin{theorem}   
Let $\Lambda$ be the set of eigenvalues of $T$, and    
let $\mathcal{B} = \{v_0,\dots,v_N\}$  be a basis of    
$\C^{N+1}$ that gives the Jordan form of $T$.    
Let $
\mathcal{H} = \bigoplus_{\lambda \in \Lambda}   
 \mathcal{H}_{\lambda} \subset \sfi,   
$ where $\mathcal{H}_{\lambda} (\varphi) = \{h \in \sfi : h \in \homog{\lambda}{k},\   
\text{for some}\ k \geq 1 \}, \lambda \in \Lambda$.   
Then we have that $\dim(\mathcal{H})= N+1$.   
\end{theorem}
\begin{remark}   
Note that we can choose both $v_0=(1,0,\dots,0)$ and $v_N=(0,\dots,0,1)$ to be in   
the basis $\mathcal{B}$, corresponding to the eigenvalues $c_0$ and   
$c_N$ respectively.   
\end{remark}   
\begin{proof}   
If $v_i \in \mathcal{B}, (0 \leq i \leq N)$,    
then $v_i \in \kker(T - \lambda I)^k$ and   
$v_i \not\in \kker(T - \lambda I)^{k-1}$, for some $\lambda \in \Lambda$,   
and $k \geq 1$. So to each $v_i \in \mathcal{B}$, we can associate   
a unique pair $(\lambda,k)$. Let us denote such $v_i = v(\lambda,k)$.   
(Note that by the previous observation, $v_0 = v(c_0,1)$   
and $v_N = v(c_N,1)$).   
   
After Theorem~\ref{tone}, we can associate to each $v(\lambda,k)$   
a function $h_{v(\lambda,k)}$ in $\homog{\lambda}{k}\cap \sfi$.   
Furthermore, the functions $\{h_{v(\lambda,k)}\}_{v \in \mathcal{B}}$,   
are linearly independent.   
For this, observe that since the vectors in $\mathcal{B}$ are   
linearly independent, its extensions $\{Y_v\}$ are   
linearly independent in $\ell(\Z)$, and therefore   
the functions $\{h_{v(\lambda,k)}\}_{v \in \mathcal{B}}$   
are linearly independent.   
   
One can see that if a finite number of functions are homogeneous   
for the same $\lambda$, then a linear combination of them   
is also homogeneous for the same $\lambda$. More precisely,   
$\sum_{i=0}^n \alpha_i h_i(\lambda,k_i) = h(\lambda, k)$,  
where  $k = \max_i(k_i)$.
  
Hence, if $p_T$, the characteristic polynomial of $T$, is factorized as:   
$p_T(x) = \prod_{\lambda \in \Lambda} (x - \lambda)^{r_\lambda}$, 
then $\dim(\mathcal{H}_{\lambda}) = r_{\lambda}$ and a basis of   
$\mathcal{H}_{\lambda}$ is the set of $(2,\lambda,k)$-homogeneous   
functions associated to the vectors $v \in \mathcal{B}$,   
such that $v = v(\lambda,k)$, for some $k \geq 1$.   
\end{proof}   

\begin{theorem} \label{main}
Assume that $\varphi$ satisfies \eqref{dilation-equation} and has
linearly independent integer translates.
Let $\mathcal{B} = \{v_0, \dots, v_N\}$ be as before, a Jordan   
basis for $T$, and let $B$ be the $(N+1)\times (N+1)$ matrix that has the   
vectors $v_i$ as rows. 
Let   
\[ \phi^0(x) = \left[\varphi(x), \varphi(x+1), \dots, 
\varphi(x+N)\right]^t \quad \text{and}\quad   
h(x) = \left[h_0(x), h_1(x), \dots,    
h_N(x)\right]^t,\]   
where $h_i$ is the homogeneous function associated to the vector $v_i$. 
Then we have
\begin{enumerate}[\upshape (i)]
\item   $h(x) = B\phi^0(x)$   $x \in [-1,1]$. \label{i}
\item   \label{ii}
$ \displaystyle
\phi^0(x) = T \phi^0(2x) \quad \text{and} 
\quad h(x) = B T B^{-1} h(2x) \quad x \in [-1/2,1/2] 
$   
where $BTB^{-1}$ is in Jordan form.  
\item \label{iii} There exists a local basis of $\sfi$ consisting of homogeneous functions. 
Moreover, if $\varphi$ has accuracy $n$, this basis can be chosen to contain the
polynomials $\{1, x, \cdots, x^{n-1}\}$.
\end{enumerate}
\end{theorem}

\begin{remark}
  Note that \eqref{ii} is a statement about the refinability of both $\phi^0$ and $h$, where
  the scaling matrix of $h$ is the Jordan form of the scaling matrix of $\phi^0$.
\end{remark}

\begin{proof} 
Since the support of $\varphi$ is $[0,N]$, $\varphi(x+k) = 0$ if 
$k \not\in \{0, 1, \dots, N\}$ for $x \in [-1,1]$. Then 
$$
h_i(x) = Y^i\Phi(x) = v_i \Phi^0(x) \ x \in [-1,1].
$$
So we have, 
\begin{equation} \label{*}
h(x) = B\phi^0(x) \ x \in [-1,1].
\end{equation}
This shows that for every interval $I \subset [-1,1]$ 
the functions $\{\varphi(x), \varphi(x + 1), \cdots, \varphi(x +N)\}$ span the same space than the functions $\{h_0, \cdots, h_N\}$ when restricted to $I$.

Since $\{\varphi(x), \varphi(x + 1), \cdots, \varphi(x +N-1)\}$ is a local
basis of \sfi, if $v_N$ has been chosen to be $v_N = (0,\cdots,0,1)$, using equation \eqref{*} $\{h_0, \cdots, h_{N-1}\}$ are a local basis for \sfi. Moreover, if $\varphi$ has accuracy $n \leq N$, 
then we can choose $v_0, \cdots, v_{n-1}$ to be the eigenvectors associated to the 
eigenvalues $\{1, \cdots, 2^{-n+1}\}$, and hence $\{h_0, \cdots, h_{N-1}\} = 
\{1, x, \cdots, x^{n-1}, h_n, \cdots, h_{N-1}\}$.

This proves \eqref{i} and \eqref{iii}.

For \eqref{ii}, using again that the support of $\varphi$ is $[0,N]$,   it is easily seen that
if $ x \in [-\frac{1}{2},\frac{1}{2}]$ then  $\phi(x) = T \phi(2x).$  Then we have
$ B\Phi^0(x) = BTB^{-1}B\Phi^0(2x) \quad x \in [-\frac{1}{2},\frac{1}{2}].
$
\end{proof}

\subsection{Generalizations} 

Throughout the paper ``function'' meant ``measurable function''. However, with the obvious modifications, all results hold for the case that $\varphi$ is a generalized function or
distribution.

It is interesting to consider the generalization of the results to arbitrary dilations $M \geq 1$ in $\R$, and in higher dimensions with arbitrary dilation matrices. It can be shown that most of the results are still true in these cases \cite{CHnM05}.   

\section{Examples for $N=3$} 
 
\noindent
{\bf B-spline} 
The simplest case of refinable functions are the {\em B-splines}. They are the 
(normalized) convolutions of the characteristic function of $[0,1]$ with itself.  
In particular the B-spline of order $3$ is the refinable function that 
satisfies:
\begin{equation} 
b(x) = \frac{1}{4} b(2x) + \frac{3}{4}b(2x-1) + \frac{3}{4} b(2x-2) 
+ \frac{1}{4} b(2x-3). 
\end{equation} 
The B-splines are those functions, for which the accuracy is maximum and so coincides with the 
dimension of the matrix $T_0$, so in this case, the eigenvalues of $T_0$ are $1$ (for the 
constant functions), $\frac{1}{2}$ (for the linear functions), 
and $\frac{1}{4}$ (for the quadratic functions). 
 
\noindent
{\bf Daubechies D$_4$} 
Daubechies wavelets, are those refinable functions of $N$ coefficients, that are 
orthogonal and provide the highest order of accuracy possible. (Note that the  
splines do not form an orthonormal basis). 
D$_4$ satisfies:
\begin{equation*} 
D_4(x) = \frac{1+\sqrt{3}}{4} D_4(2x) + \frac{3+\sqrt{3}}{4}D_4(2x-1) + 
\frac{3-\sqrt{3}}{4} D_4(2x-2) + \frac{1-\sqrt{3}}{4} D_4(2x-3). 
\end{equation*} 
D$_4$ has accuracy 2 (it reproduces the constant and the linear functions). 
In this case the matrix $T_0$ has eigenvalues $1$, $\frac{1}{2}$ and 
$c_0 = \frac{1+\sqrt{3}}{4}$. So a basis for $\text{span}\{D_4(x),%
 D_4(x-1), D_4(x-2)\}_{x \in [0,1]}$ is also given by 
$\text{span}\{1, x, h_{c_0}(x)\}_{x \in [0,1]}$ where $h_{c_0}$ is the homogeneous function associated to $c_0$. 
 
\noindent
{\bf $(\lambda,1)$-Homogeneous functions are not enough} 
In the two previous examples, we could always obtain a basis of  
$\text{span}\{f(x), f(x-1), f(x-2)\}_{x \in [0,1]}$ just by using 1-homogeneous functions. The following example is to illustrate, that even in the simple case of only 4 coefficients, 
it may be necessary to use homogeneous functions of order bigger than 1. 
Consider the function: 
\begin{equation} 
f(x) = \frac{1}{3} f(2x) + \frac{2}{3}f(2x-1) + \frac{2}{3} f(2x-2) 
+ \frac{1}{3} f(2x-3). 
\end{equation} 
It can be shown that $f$ has accuracy 1, and the eigenvalues of $T$ are $\{1,\frac{1}{3}\}$. 
So in this case, $\text{span}\{f(x), f(x-1), f(x-2)\}_{x \in [0,1]} = 
\text{span}\{1, h_{\{1/3,1\}}(x), h_{\{1/3,2}\}(x)\}_{x \in [0,1]}$, where $h_{\{1/3,1\}}$ is a 1-homogeneous function corresponding to the eigenvalue $1/3$, and $h_{\{1/3,2\}}$ is a 2-homogeneous function corresponding to the eigenvalue $1/3$. 
 
 \section{Acknowledgments} We wish to thank the anonymous referee for many helpful suggestions.
 
   
\bibliographystyle{amsalpha}   
\bibliography{cyu}   
   
\end{document}